\documentclass{amsart}
\usepackage[bookmarksnumbered,colorlinks, linkcolor=blue, citecolor=red, pagebackref, bookmarks, breaklinks]{hyperref}

\newtheorem{thm}{Theorem}[section]
\newtheorem{cor}[thm]{Corollary}
\newtheorem{lema}[thm]{Lemma}
\newtheorem{prop}[thm]{Proposition}
\theoremstyle{definition}
\newtheorem{defn}[thm]{Definition}

\theoremstyle{remark}
\newtheorem{rem}[thm]{Remark}

\numberwithin{equation}{section}
\newcommand{\R}{\mathbb R}
\newcommand{\N}{\mathbb{N}}
\def\diver{\mathop{\text{\normalfont div}}}

\begin{document}

\title[Eigenvalues of the fractional $g$-Laplacian]{Variational Eigenvalues of the fractional $g$-Laplacian}

\author[S.Bahrouni]{Sabri Bahrouni}
\author[H.Ounaies]{Hichem Ounaies}
\author[A.Salort]{Ariel Salort}
\address[S.Bahrouni and H.Ounaies]{Mathematics Department, Faculty of Sciences, University of Monastir, 5019 Monastir, Tunisia}
\address[A. Salort]{IMAS - CONICET, Ciudad Universitaria, Pabell\'on I (1428) Av. Cantilo s/n, Buenos Aires, Argentina}
\email[S.Bahrouni]{sabri.bahrouni@fsm.rnu.tn}
\email[H.Ounaies]{hichem.ounaies@fsm.rnu.tn}
\email[A. Salort]{asalort@dm.uba.ar}

\keywords{Fractional Orlicz-Sobolev spaces, Ljusternik-Schnirelman eigenvalues, Minimax eigenvalues \\
\hspace*{.3cm} {\it 2010 Mathematics Subject Classifications}:
 46E30, 35R11, 45G05}

\begin{abstract}
In the present work we study existence of sequences of variational eigenvalues to non-local non-standard growth  problems ruled by the fractional $g-$Laplacian operator with different boundary conditions (Dirichlet, Neumann and Robin). Due to the non-homogeneous nature of the operator several drawbacks must be overcome, leading to some results that contrast with the case of power functions.
\end{abstract}

\maketitle
\tableofcontents

\section{Introduction}

Given two functionals $\mathcal{A}$ and $\mathcal{B}$ defined on a suitable space $\mathcal{X}$ and a prescribed number $c$, the task of analyzing the existence of numbers $\lambda\in\R$ and elements $u\in \mathcal{X}$ satisfying (in some appropriated sense) equations of the type
$$
\mathcal{B}'(u)=\lambda \mathcal{A}'(u), \qquad \mathcal{A}(u)=c,
$$
has been a challenging labor whose beginning  dates back to the mid-20th century (here $\mathcal{A}'$ and $\mathcal{B}'$ denote the Fr\'echet derivatives of the functionals). The study on Hilbert spaces was addressed  by Krasnoselskij in \cite{K64}; for Banach spaces, it can be found in the literature the works of Citlanadze \cite{C53} and Browder \cite{B65, Br53, Br65b}, where the notion of \emph{category of sets} in the sense of  Ljusternik and Schnirelman is used. See also \cite{FN72,FN73} for some applications to partial differential equations. The amount of research on these topics is nowadays huge. For practical reasons, for an introduction to this theory and a comprehensive list of references we refer to   the books \cite{B68, Br65, Papageorgiou, zeidlerIII}.

After the introduction of the so-called monotonicity methods by Browder \cite{B63,B64}, Minty \cite{M62} and Va\v{i}nberg and Ka\v{c}urovski\v{i} \cite{VK59}, the study of quasilinear operators experimented an explosive growth, and both variational and nonvariational techniques were introduced by Browder, Fu\v{c}\'{\i}k, Ladyzhenskaya, Leray, Lions, Morrey, Ne\v{c}as, Rabinowicz, Schauder, Serrin, and Trudinger, among several other mathematicians.

The prototypical $p-$Laplace operator ($p>1$) then became a focus of study, and in particular, to understand its spectral structure:  given an open and bounded set $\Omega\subset \R^n$,  to determine the existence  couples $(\lambda,u)$ satisfying the equation
\begin{equation} \label{p.lap}
-\diver(|\nabla u|^{p-2}u)=\lambda |u|^{p-2}u \text{ in } \Omega, \quad u=0 \text{ on } \partial \Omega
\end{equation}
in a suitable sense. In the seminal work  of Garc\'{\i}a Azorero and Peral Alonso \cite{GP}, it was  proved the existence of a variational sequence of eigenvalues tending to $+\infty$, which however, it is not known to exhaust the spectrum unless $p=2$ or $n=1$. Several properties on eigenvalues (and their corresponding eigenfunctions) were addressed by Anane et al \cite{A87, A94} and Lindqvist \cite{L90}, among others, and also for more general boundary conditions than  Dirichlet. See also \cite{FBPS, An06}.

At this point,  two possible generalizations of the eigenvalue problem \eqref{p.lap} could be considered. First, its non-local counterpart governed by the well-known \emph{fractional $p-$Laplace} operator  takes the form
\begin{equation} \label{frac.p.lap}
(-\Delta_p)^s u =\lambda |u|^{p-2}u \text{ in } \Omega, \quad u=0 \text{ in } \R^n \setminus \Omega
\end{equation}
where $s\in(0,1)$ is a fractional parameter, $p>1$, and
$$
(-\Delta_p)^s u(x):=\text{p.v.} \int_{\R^n} \frac{|u(x)-u(y)|^{p-2}(u(x)-u(y))}{|x-y|^{n+sp}}\,dy.
$$
The main difference here arises in the fact that this operator takes into account   interactions coming from the whole space. The same occurs with the \emph{boundary} condition. Problem \eqref{frac.p.lap} was introduced in \cite{LL14}. Existence of a sequence of (variational) eigenvalues to \eqref{frac.p.lap} and its behavior as $s\uparrow 1$ was dealt in \cite{BPS16}. Several properties on eigenvalues and eigenfunctions were obtained in \cite{BP16, DPS15, FP14, LL14}.

A second possible generalization of \eqref{p.lap} can be obtained keeping the local structure of the operator but allowing a growth behavior more general than a power. These considerations lead to the well-known \emph{$g-$Laplace} operator and the problem
$$
-\Delta_g u = \lambda g(|u|)\frac{u}{|u|} \text{ in } \Omega, \quad u=0 \text{ on } \partial\Omega,
$$
where $\Delta_g u:= \diver\left(g(|\nabla u|)\frac{\nabla u}{|\nabla u|}\right)$. The function $g=G'$ is given in terms of a so-called Young function $G$. Here the structure of the spectrum radically changes due to the non-homogeneity of the operator, fact that is key in the arguments corresponding to results related to problems \eqref{p.lap} and \eqref{frac.p.lap}. Here several important differences appear: the spectrum may not be discrete, then it is not clear the notion of a first eigenvalue nor of a sequence of them. However, when restricting the \emph{energy level} of functions, some properties of operators with $p-$structure are recovered: in \cite{GM02,GHLMS99,MT99, ML12} existence of eigenvalues was studied, and in \cite{E96, T00} the existence of a discrete sequence of eigenvalues was obtained. It is worth of mention that under these settings, eigenvalues are not in general variational, so less information can be retrieved when using standard techniques.

The object of study of this manuscript is the eigenvalue problem ruled by the \emph{fractional $g-$Laplacian}, which can be seen as the non-local non-standard growth counterpart of \eqref{p.lap}. Given an open and bounded set  $\Omega\subset\R^n$ and a fractional parameter $s\in (0,1)$, we are concerned with the following eigenvalue problems:
\begin{itemize}
\item Dirichlet problem:
\begin{align}\tag{$D(\Omega)$} \label{eq.d}
\begin{cases}
(-\Delta_g)^s u =\lambda g(|u|)\frac{u}{|u|} &\text{ in } \Omega\\
u=0 &\text{ in } \R^n \setminus \Omega,
\end{cases}
\end{align}

\item Neumann problem:
\begin{align}\tag{$N(\Omega)$} \label{eq.n}
\begin{cases}
(-\Delta_{g})^s u =\lambda g(|u|)\frac{u}{|u|} &\text{ in } \Omega\\
\mathcal{N}_g u =0 &\text{ in } \R^n \setminus \Omega.
\end{cases}
\end{align}

\item Robin problem: given $\beta\in L^\infty(\R^n\setminus \Omega)$
\begin{align}\tag{$R(\Omega)$} \label{eq.r}
\begin{cases}
(-\Delta_{g})^s u =\lambda g(|u|)\frac{u}{|u|} &\text{ in } \Omega\\
\mathcal{N}_g u + \beta g(|u|)\frac{u}{|u|}=0 &\text{ in } \R^n \setminus \Omega.
\end{cases}
\end{align}
\end{itemize}
The precise notion of \emph{eigenvalues} and \emph{eigenfunctions} to these problems is defined in Section \ref{section.lagrange}.

Here, the fractional $g-$Laplacian is defined as
\begin{equation}
(-\Delta_g)^s u:= \, \text{p.v.} \int_{\mathbb{R}^n} g\left( |D_s u|\right)\frac{D_s u}{|D_s u|} \frac{dy}{|x-y|^{n+s}},
\end{equation}
where $G$ is a Young function (see Section \ref{sec.prel} for the precise definition) such that $g=G'$ and $D_s u:=\frac{u(x)-u(y)}{|x-y|^s}$.  Clearly, $(-\Delta_{g})^s$ boils down to the fractional $p-$Laplacian when $G$ is a power function and to the $p-$Laplacian when $s\uparrow 1$. See \cite{FBS}. The boundary conditions in the previous problems   reflect the non-local nature of the operator $(-\Delta_{g})^s$: the Dirichlet case corresponds to functions vanishing outside $\Omega$ and not only on $\partial\Omega$, whereas the Neumann and Robin equations  make use of the nonlocal normal derivative $\mathcal{N}_g$ introduced in \cite{BS20}.

Throughout this article we assume the Young function $G=\int_0^tg(t)\,dt$   to satisfy the following structural conditions:
\begin{equation} \label{G1} \tag{$G_1$}
1<p^-\leq \frac{tg(t)}{G(t)} \leq p^+<\infty \quad \forall t>0,
\end{equation}
\begin{equation} \label{G2} \tag{$G_2$}
t\mapsto G(\sqrt{t}),\ t\in[0,\infty)  \text{ is convex},
\end{equation}
\begin{equation} \label{G3} \tag{$G_3$}
 \displaystyle\int_{0}^{1}\frac{G^{-1}(\tau)}{\tau^{\frac{n+s}{n}}}d\tau<\infty\ \quad \text{and}\quad
\displaystyle\int_{1}^{+\infty}\frac{G^{-1}(\tau)}{\tau^{\frac{n+s}{n}}}d\tau=\infty.
\end{equation}
In the case of powers (i.e. $G(t)=t^p$) these conditions mean that $2\leq p \leq \frac{np}{n-sp}$.

The following functionals have an essential role in the study of eigenvalue problems for the fractional $g-$Laplacian
\begin{equation} \label{func.1}
\mathcal{I}(u):=\int_\Omega G(|u|)\,dx, \qquad \mathcal{J}(u):=\iint_{\R^n\times\R^n} G(|D_s u|)\frac{dxdy}{|x-y|^n}.
\end{equation}
and
\begin{equation} \label{func.2}
\mathcal{J}_\beta(u):=\iint_{\R^{2n}\setminus (\Omega^c)^2} G(|D_s u|)\frac{dxdy}{|x-y|^n} + \int_{\R^n \setminus \Omega} \beta G(|u|)\,dx
\end{equation}
(see Section \ref{section.lagrange} for details).

Problem \ref{eq.d} was recently studied in \cite{VS20,S19}, where, once established  existence of a minimizer of the associated problem under energy constraint
$$
\Lambda_\alpha^D:=\min\left\{ \mathcal{J}(u) \colon \mathcal{I}(u) =\alpha\right\}
$$
for a given value $\alpha>0$, then by means of a version of the Lagrange multipliers rule  it is deduced existence of an eigenvalue $\lambda_\alpha^D$ of \ref{eq.d} in a suitable sense, and in general $\Lambda_\alpha^D\neq \lambda_\alpha^D$. A further extension of this result for problems \ref{eq.n} and \ref{eq.r} is provided in \cite{BS20}. Other very recent results involving eigenvalues of the fractional $g-$Laplacian can be found in \cite{Cianchi, Cianchi2, ABS, Sabri1, Sabri2, Sabri3, Sabri5,     DNFBS,  FBPLS,FBSV, VS20, S20}.


In the wake of the non-homogeneity of the operator, many of the properties that eigenvalues of the fractional $p-$Laplacian fulfill are possibly not  inherited in the non-homogeneous case: for instance, isolation, simplicity and a variational characterization of the first eigenvalue, or a variational formula for the second one (see \cite{DPS15, FP14, LL14}). Moreover,  the spectrum of \ref{eq.d} could be continuous when $G$ is a general Young function, and in principle, it is not clear the meaning of a first or second eigenvalue. Due to these drawbacks, the main aim of this manuscript is to understand under which conditions it is possible to build a sequence of eigenvalues to problems \ref{eq.d}, \ref{eq.n} and \ref{eq.r}.

Our arguments are based in the fact that existence of sequences of variational eigenvalues can be established when prescribing some energy level. First, when the quantities $\mathcal{J}(u)$ or $\mathcal{J}_\beta(u)$ involving  the $s-$H\"older quotient  are prescribed, by means of the Ljusternik-Schnirelman theory we can infer  existence of a discrete sequence of critical points of $\mathcal{I}(u)$ which allow us to build a sequence of non-negative eigenvalues. With the notation introduced in Section \ref{sec.prel}, our first result reads as follows.

\begin{thm}\label{main.result.1}
For any $\alpha>0$ there exist sequences of non-negative   numbers $\{\lambda^D_{k,\alpha}\}_{k\in\mathbb{N}}, \{\lambda^N_{k,\alpha}\}_{k\in\mathbb{N}}$ and $\{\lambda^R_{k,\alpha}\}_{k\in\mathbb{N}}$ which are eigenvalues of problems \ref{eq.d}, \ref{eq.n} and \ref{eq.r}, respectively. Moreover, these sequences diverge as $k\to\infty$.

The corresponding eigenfunctions
$$
\{u_{k,\alpha}^D\}_{k\in \mathbb{N}}\subset W^{s,G}_0(\Omega), \quad \{u_{k,\alpha}^N\}_{k\in \mathbb{N}}\subset W^{s,G}_*(\Omega) \quad \text{ and } \quad
\{u_{k,\alpha}^R\}_{k\in \mathbb{N}}\subset \mathcal{X}_\beta(\Omega)
$$
satisfy the constraints
$$
\mathcal{J}(u_{k,\alpha}^D)=\mathcal{J}_0(u_{k,\alpha}^N)=\mathcal{J}_\beta(u_{k,\alpha}^R)=\alpha
$$
and
$$
\mathcal{I}(u_{k,\alpha}^D)=c_{k,\alpha}^D, \qquad \mathcal{I}(u_{k,\alpha}^N)=c_{k,\alpha}^N, \qquad  \mathcal{I}(u_{k,\alpha}^R)=c_{k,\alpha}^R
$$
where the critical values are obtained as
$$
c_{k,\alpha}^D = \sup_{K\in \mathcal{C}_k^D}\inf_{u\in K}\mathcal{I}(u), \qquad
c_{k,\alpha}^N = \sup_{K\in \mathcal{C}_k^N}\inf_{u\in K}\mathcal{I}(u), \qquad
c_{k,\alpha}^R = \sup_{K\in \mathcal{C}_k^R}\inf_{u\in K}\mathcal{I}(u)
$$
where, denoting by $\gamma(K)$ the Krasnoselskii genus of $K$, we define
$$
\mathcal{C}_k^D:=\{ K\subset M_\alpha^D \text{ compact, symmetric with } \mathcal{I}(u)>0 \text{ on } K  \text{ and } \gamma(K)\geq k\},
$$
and
\begin{align*}
M_\alpha^D &:= \{ u\in W^{s,G}_0(\Omega) \colon \mathcal{J}(u)=\alpha\},\\
M_\alpha^N &:= \{ u\in W^{s,G}_*(\Omega) \colon \mathcal{J}_0(u)=\alpha\},\\
M_\alpha^R &:= \{ u\in \mathcal{X}_\beta(\Omega) \colon \mathcal{J}_\beta (u)=\alpha\}.
\end{align*}
The sets $\mathcal{C}_k^N$ and $\mathcal{C}_k^R$ are defined analogously by changing the superscript $D$ by $N$ and $R$, respectively.
\end{thm}

Secondly, when  the quantity $\mathcal{I}(u)$ is prescribed, by using the  so-called minimax theory we obtain a sequence of critical points of $\mathcal{J}(u)$ which provides for a (different) sequence of eigenvalues:

\begin{thm}\label{main.result.2}
For any $\alpha>0$ there exist  sequences of non-negative   eigenvalues $\{\Lambda^D_{k,\alpha}\}_{k\in\mathbb{N}}$, $\{\Lambda^N_{k,\alpha}\}_{k\in\mathbb{N}}$ and $\{\Lambda^R_{k,\alpha}\}_{k\in\mathbb{N}}$ of \ref{eq.d}, \ref{eq.n} and \ref{eq.r}, respectively.

The corresponding eigenfunctions
$$
\{u_{k,\alpha}^D\}_{k\in \mathbb{N}}\subset W^{s,G}_0(\Omega), \quad \{u_{k,\alpha}^N\}_{k\in \mathbb{N}}\subset W^{s,G}_*(\Omega) \quad \text{ and } \quad
\{u_{k,\alpha}^R\}_{k\in \mathbb{N}}\subset \mathcal{X}_\beta(\Omega)
$$
satisfy the constraints
$$
\mathcal{I}(u_{k,\alpha}^D)=\mathcal{I}(u_{k,\alpha}^N)=\mathcal{I}(u_{k,\alpha}^R)=\alpha
$$
and
$$
\mathcal{J}(u_{k,\alpha}^D)=C_{k,\alpha}^D, \qquad \mathcal{J}_0(u_{k,\alpha}^N)=C_{k,\alpha}^N, \qquad  \mathcal{J}_\beta (u_{k,\alpha}^R)=C_{k,\alpha}^R
$$
where the critical values are obtained as
$$
C_{k,\alpha}^D = \inf_{h\in \Gamma(S^{k-1},M_\alpha^D)} \sup_{w \in S^{k-1}} \mathcal{J}(h(w))
$$
and $\Gamma(S^k,M_\alpha^D)=\{h \in C(S^k,M_\alpha^D)\colon h \text{ is odd}\}$, being $S^k$ the unit sphere in $\R^{k+1}$.

The numbers $C_{k,\alpha}^N$ and $C_{k,\alpha}^R$ are defined analogously by changing the superscript $D$ by $N$ and $R$, respectively, and the functional $\mathcal{J}$ by  $\mathcal{J}_0$ and $\mathcal{J}_\beta$, respectively.  $M_\alpha^D$, $M_\alpha^N$ and $M_\alpha^R$ denote  the sets
\begin{align*}
M_\alpha^D &:= \{ u\in W^{s,G}_0(\Omega) \colon \mathcal{I}(u)=\alpha\},\\
M_\alpha^N &:= \{ u\in W^{s,G}_*(\Omega) \colon \mathcal{I}(u)=\alpha\},\\
M_\alpha^R &:= \{ u\in \mathcal{X}_\beta(\Omega) \colon \mathcal{I}(u)=\alpha\}.
\end{align*}
\end{thm}


In contrast to what happens in the case of powers,  the Ljusternik-Schnirelman and minimax  eigenvalues obtained in Theorems \ref{main.result.1} and \ref{main.result.2}, in general, are neither the same nor easily  comparable each other. As a consequence, in principle we cannot ensure that the minimax sequences built in Theorems \ref{main.result.2} diverge as $k\to\infty$. For a discussion of this situation in the case $G(t)=t^p$ see \cite{DR, LS}. However, in the Dirichlet case we are able to compare miminax eigenvalues of the $g-$Laplacian with eigenvalues of the $p^--$Laplacian (where $p^-$ is given in \eqref{G1}), giving as a consequence the following:

\begin{thm} \label{dirichlet.diverge}
With the notation of Theorem \ref{main.result.2},
$$
\Lambda_{k,\alpha}^D\to\infty \text{ as }k\to\infty, \qquad
C_{k,\alpha}^D\to\infty \text{ as }k\to\infty.
$$
\end{thm}

In the case of powers,  eigenvalues and critical points obtained in Theorems \ref{main.result.1} and \ref{main.result.2} coincide, i.e.,   $\lambda_k=c_k$ and $\Lambda_k=C_k$ for all $k\in\N$. It does not occur for a general Young function, although the following comparison result holds.

\begin{thm} \label{compara}
Let $\{\lambda_{k,\alpha}^D\}_{k\in\N}$ and $\{c_{k,\alpha}^D\}_{k\in\N}$ be as in Theorem \ref{main.result.1}, and let $\{\Lambda_{k,\alpha}^D\}_{k\in\N}$ and $\{C_{k,\alpha}^D\}_{k\in\N}$ be as in Theorem \ref{main.result.2}. Then
$$
\frac{ \alpha p^-}{p^+}  \leq c_{k,\alpha}^D \cdot \lambda_{k,\alpha}^D \leq
\frac{ \alpha p^+}{p^-}  , \qquad
\frac{p^-}{ \alpha p^+}  \leq \frac{\Lambda_{k,\alpha}^D}{C_{k,\alpha}^D} \leq
\frac{p^+}{ \alpha p^-},
$$
where $p^+$ and $p^-$ are the numbers defined in \eqref{G1}.
\end{thm}

Finally, our last result establishes the closedness of the spectrum of the fractional $g-$Laplacian in the following sense:
\begin{thm} \label{lam.inf}
Fixed $\alpha_0>0$, let $\{\alpha_k\}_{k\in\N}\subset (0,\alpha_0)$. Let $\{\lambda_{\alpha_k}\}_{k\in\N}$ be a sequence of eigenvalues of  \ref{eq.d} with eigenfunctions $\{u_{\alpha_k}\}_{k\in\N}\subset W^{s,G}_0(\Omega)$ such that $\Phi_{G,\Omega}(u_{\alpha_k})=\alpha_k$.

Then, if $\lim_{k\to\infty} \lambda_{\alpha_k}=\lambda$, we have that $\lambda$ is an eigenvalue of \ref{eq.d} with eigenfunction $u\in W^{s,G}_0(\Omega)$ such that $\Phi_{G,\Omega}(u)=\beta\in (0,\alpha_0)$.

A similar assertion holds for sequences of eigenvalues of \ref{eq.n} and \ref{eq.r}.
\end{thm}

As a consequence we get the following.
\begin{cor}
Fixed $\alpha_0>0$, the numbers
$$
\hat \lambda_k^D:=\inf_{\alpha\in (0,\alpha_0)} \lambda_{k,\alpha}^D, \quad
\hat \lambda_k^N:=\inf_{\alpha\in (0,\alpha_0)} \lambda_{k,\alpha}^N,\quad
\hat \lambda_k^R:=\inf_{\alpha\in (0,\alpha_0)} \lambda_{k,\alpha}^R
$$
are eigenvalues of \ref{eq.d}, \ref{eq.n} and \ref{eq.r}, respectively.

The same claim is true for the numbers $\hat \Lambda_{k}^D$, $\hat  \Lambda_{k}^N$ and $\hat  \Lambda_{k}^R$ defined in an analogous way.
\end{cor}

\section{Definitions and preliminary results} \label{sec.prel}

In this section we introduce the classes of Young function and fractional Orlicz-Sobolev functions, the
suitable class where the fractional $g$-Laplacian is well defined.

\subsection{Young functions}

An application $G\colon\R_+\to \R_+$ is said to be a  \emph{Young function} if it admits the integral formulation $G(t)=\int_0^t g(\tau)\,d\tau$, where the right continuous function $g$ defined on $[0,\infty)$ has the following properties:
\begin{align*}
&g(0)=0, \quad g(t)>0 \text{ for } t>0 \label{g0} \tag{$g_1$}, \\
&g \text{ is non-decreasing on } (0,\infty) \label{g2} \tag{$g_2$}, \\
&\lim_{t\to\infty}g(t)=\infty  \label{g3} \tag{$g_3$} .
\end{align*}
From these properties it is easy to see that a Young function $G$ is continuous, nonnegative, strictly increasing and convex on $[0,\infty)$.

The following  properties on Young functions are well-known. See for instance \cite{FJK} for the proof of these results.

\begin{lema} \label{lema.prop}
Let $G$ be a Young function satisfying \eqref{G1} and $a,b\geq 0$. Then
\begin{align*}
  &\min\{ a^{p^-}, a^{p^+}\} G(b) \leq G(ab)\leq   \max\{a^{p^-},a^{p^+}\} G(b),\tag{$L_1$}\label{L1}\\
  &G(a+b)\leq \mathbf{C} (G(a)+G(b)) \quad \text{with } \mathbf{C}:=  2^{p^+},\tag{$L_2$}\label{L2}\\
	&G \text{ is Lipschitz continuous}. \tag{$L_3$}\label{L_3}
 \end{align*}
\end{lema}
Condition \eqref{G1} is known as the \emph{$\Delta_2$ condition} or \emph{doubling condition} and, as it is showed in \cite[Theorem 3.4.4]{FJK}, it is equivalent to the right hand side inequality in \eqref{G1}.

The \emph{complementary Young function} $\tilde G$ of a Young function $G$ is defined as
$$
\tilde G(t):=\sup\{tw -G(w): w>0\}.
$$
From this definition the following Young-type inequality holds
\begin{equation} \label{Young}
ab\leq G(a)+\tilde G(b)\qquad \text{for all }a,b\geq 0,
\end{equation}
and the following H\"older's type inequality
$$
\int_\Omega |uv|\,dx \leq \|u\|_{L^G(\Omega)} \|v\|_{L^{\tilde G}(\Omega)}
$$
for all $u\in L^G(\Omega)$ and $v\in L^{\tilde G}(\Omega)$. Moreover, it is not hard to see that $\tilde G$ can be written in terms of the inverse of $g$ as
\begin{equation} \label{xxxx}
\tilde G(t)=\int_0^t g^{-1}(\tau)\,d\tau,
\end{equation}
see \cite[Theorem 2.6.8]{RR}.

Since $g^{-1}$ is increasing, from \eqref{xxxx} and \eqref{G1} it is immediate the following relation.
\begin{lema} \label{lemita.1}
Let $G$ be an Young function satisfying \eqref{G1} such that $g=G'$ and denote by  $\tilde G$ its complementary function. Then
$$
\tilde G(g(t)) \leq p^+ G(t)
$$
holds for any $t\geq 0$.
\end{lema}

The following convexity property proved in \cite{Lamperti}[Lemma 2.1] will be useful.
\begin{lema} \label{lemita}
Let $G$ be a Young function satisfying \eqref{G1} and \eqref{G2}. Then for every $a,b\in \R$,
$$
\frac{G(|a|) + G(|b|)}{2} \geq
    G\left(\left|\frac{a+b}{2} \right| \right) + G\left(\left|\frac{a-b}{2} \right| \right).
$$
\end{lema}

\subsection{Fractional Orlicz-Sobolev spaces}

Given a Young function $G$ such that $G'=g$, a parameter $s\in(0,1)$ and an open and bounded set $\Omega\subseteq \R^n$ we consider the spaces
\begin{align*}
&L^G(\Omega) :=\left\{ u\colon \Omega \to \R \text{ measurable  such that }  \Phi_{G,\Omega}(u) < \infty \right\},\\
&W^{s,G}(\R^n):=\left\{ u\in L^G(\R^n) \text{ such that } \Phi_{s,G,\R^n}(u)<\infty \right\},\\
&W^{s,G}_0(\Omega):=\left\{ u\in W^{s,G}(\R^n) :\ u=0\ a.e.\ \text{in}\ \R^n\setminus\Omega\right\}
\end{align*}
where the modulars $\Phi_{G,\Omega}$ and $\Phi_{s,G,\R^n}$ are defined as
\begin{align*}
&\Phi_{G,\Omega}(u):=\int_{\Omega} G(|u(x)|)\,dx\\
&\Phi_{s,G,\R^n}(u):=
  \iint_{\R^n\times\R^n} G( |D_su(x,y)|)  \,d\mu,
\end{align*}
where  the \emph{$s-$H\"older quotient} is defined as
$$
D_s u(x,y):=\frac{u(x)-u(y)}{|x-y|^s},
$$
being $d\mu(x,y):=\frac{ dx\,dy}{|x-y|^n}$.
These spaces are endowed with the so-called \emph{Luxemburg norms}
\begin{align*}
&\|u\|_{L^G(\Omega)} := \inf\left\{\lambda>0\colon \Phi_{G,\Omega}\left(\frac{u}{\lambda}\right)\le 1\right\},\\
&\|u\|_{W^{s,G}(\Omega)} := \|u\|_{L^G(\Omega)} + [u]_{W^{s,G}(\R^n)},
\end{align*}
where the  {\em $(s,G)$-Gagliardo semi-norm} is defined as
\begin{align*}
&[u]_{W^{s,G}(\R^n)} :=\inf\left\{\lambda>0\colon \Phi_{s,G,\R^n}\left(\frac{u}{\lambda}\right)\le 1\right\}.
\end{align*}
Since we assume \eqref{G1}, the space $W^{s,G}(\Omega)$ is a reflexive Banach space. Moreover $C_c^\infty$ is dense in $W^{s,G}(\R^n)$. See \cite[Proposition 2.11]{FBS} and \cite[Proposition 2.9]{DNFBS} for details.

Observe that, in light of the modular Poincar\'e inequality for $W^{s,G}(\Omega)$ (see \cite{S19}), it follows that $[\,.\,]_{W^{s,G}(\R^n)}$ is an equivalent norm in $W_0^{s,G}(\Omega)$. We observe that $W^{s,G}_0(\Omega)$ is the natural space to deal with Dirichlet problems.

The space of fractional Orlicz-Sobolev functions is the appropriated one where to define the \emph{fractional $g-$Laplacian operator}
$$
(-\Delta_g)^s u :=2 \,\text{p.v.} \int_{\R^n} g( |D_s u|) \frac{D_s u}{|D_s u|} \frac{dy}{|x-y|^{n+s}},
$$
where \text{p.v.} stands for {\em in principal value}. This operator  is well defined between $W^{s,G}(\R^n)$ and its dual space $W^{-s,\tilde G}(\R^n)$ (see \cite{FBS} for details). In fact, it follows that
$$
\langle (-\Delta_g)^s u,v \rangle =   \iint_{\R^n\times\R^n} g(|D_s u|) \frac{D_s u}{|D_s u|}  D_s v \,d\mu,
$$
for any $v\in W^{s,G}(\R^n)$.

As proved in \cite[Proposition 2.6]{BS20}, the following integration by parts formula arise naturally for $u\in C^2$ functions
$$
\langle (-\Delta_g)^s u,v \rangle_{*} = \int_{\Omega} v  (-\Delta_g)^s u \ dx + \int_{\mathbb{R}^n \setminus \Omega} v  \mathcal{N}_g u \ dx  \quad  \forall v\in C^2
$$
where the \emph{normal derivative} $\mathcal{N}_g$ is defined as
$$
\mathcal{N}_g u(x) = \int_\Omega g(|D_s u|) \frac{D_s u}{|D_s u|}\frac{dy}{|x-y|^{n+s}}.
$$
and the product $\langle \cdot ,\cdot \rangle_*$ is defined as
$$
\langle (-\Delta_g)^s u,v \rangle_* =    \iint_{\mathbb{R}^{2n} \setminus (\R^n\setminus \Omega)^2} g(|D_s u|) \frac{D_s u}{|D_s u|}  D_s v \,d\mu.
$$
The previous definitions induce the following notation
\begin{align*}
\Phi_{s,G,*}(u)&=\iint_{\R^{2n}\setminus (\Omega)^2} G(|D_s u(x,y)|)\,d\mu\\
[u]_{W^{s,G}_*(\R^n)} &= \inf\left\{ \lambda>0\colon  \Phi_{s,G,*}\left(\frac{u}{\lambda}\right) \leq 1 \right\}.
\end{align*}
Hence, it is natural to defined the space
$$
W^{s,G}_*(\Omega):=\{u\in L^G(\Omega)   \text{ such that } \Phi_{s,G,*}(u)<\infty\},
$$
which will be the appropriated one when dealing with the Neumann boundary condition.

The suitable space in which to define a Robin boundary condition of the type $\mathcal{N}_g u + \beta g(|u|)u/|u|$ in $\R^n \setminus \Omega$ (where $\beta$ is a fixed function in $L^\infty(\R^n \setminus \Omega)$) is:
$$
\mathcal{X}_\beta(\Omega)=\left\{ u \text{ measurable} \colon \Phi_{s,G,*}(u) + \Phi_{G,\Omega}(u)+ \Phi_{G, \beta,  \R^n\setminus \Omega}(u) <\infty \right\}
$$
where
$$
\Phi_{G,\beta,\R^n\setminus\Omega}(u) = \int_{\R^n\setminus\Omega} \beta G(|u|)\,dx.
$$
This space can be proved to be a reflexive Banach space endowed with the norm
$$
\|u\|_\mathcal{X_\beta}:=[u]_{W^{s,G}_*(\R^n)}+ \|u\|_{L^G(\Omega)} +\| u\|_{L^{G,\beta}(\R^n\setminus \Omega)},
$$
being
$$
\| u\|_{L^{G,\beta}(\R^n\setminus \Omega)}= \inf\left\{ \lambda>0\colon \int_{\mathbb{R}^n\setminus\Omega} \beta\ G\left(\frac{|u|}{\lambda}\right)\,d\mu \leq 1 \right\}.
$$
See \cite{BS20} for more details.

\begin{rem}
Observe that, with the notation introduced in this section, the functionals introduced in \eqref{func.1} and \eqref{func.2} can be identified as
\begin{align} \label{funcionales}
\begin{split}
&\mathcal{I}\colon W^{s,G}(\Omega)\to \R, \quad \mathcal{I}(u)=\Phi_{G,\Omega}(u),\\
&\mathcal{J}\colon W^{s,G}_0(\Omega)\to \R, \quad \mathcal{J}(u)=\Phi_{s,G,\R^n}(u),\\
&\mathcal{J}_\beta \colon (\Omega)\to \R, \quad \mathcal{J}_\beta(u)=\Phi_{s,G,*}(u) + \Phi_{G,\beta,\R^n\setminus \Omega}.
\end{split}
\end{align}
\end{rem}

In order to state some embedding results for fractional Orlicz-Sobolev spaces we recall that given two Young functions $A$ and $B$, we say that \emph{$B$ is essentially stronger than $A$} or equivalently that \emph{$A$ decreases essentially more rapidly than $B$}, and denoted by $A\prec \prec B$, if for each $a>0$ there exists $x_a\geq 0$ such that $A(x)\leq B(ax)$ for $x\geq x_a$.

When the Young function $G$ fulfills condition \eqref{G3}, the critical function for the fractional Orlicz-Sobolev embedding is given by
$$
G_{*}^{-1}(t)=\int_{0}^{t}\frac{G^{-1}(\tau)}{\tau^{\frac{n+s}{n}}}d\tau.
$$

With these preliminaries the following compact embeddings hold.

\begin{prop}\label{ceb}
	Let $G$ be a Young function satisfying \eqref{G1} and \eqref{G3} and let $s\in(0,1)$. Let $\Omega\subset \R^n$ be a $C^{0,1}$ bounded open subset. Then for any Young function $B$ such that $B \prec \prec G_{*}$ it holds that
\begin{itemize}
  \item[(i)]
  the embedding  $W^{s,G}(\Omega)\hookrightarrow L^{B}(\Omega)$ is compact;

  \item[(ii)]
  the embedding $\mathcal{X}_\beta(\Omega) \hookrightarrow L^{B}(\Omega)$ is compact.
\end{itemize}
\end{prop}
The proof of (i) can be found in \cite{Cianchi}[Theorem 9.1], \cite{Sabri2}[Theorem 1.2]; for the proof of (ii) see \cite{BS20}[Lemma 3.2].

The following relation between modulars and norms holds. See  \cite[Lemma 2.1]{Fukagai}.
\begin{lema}\label{ineq1}
   Let $G$ be a Young function satisfying \eqref{G1} and let $\xi^-(t)=\min\{t^{p^-},t^{p^+}\}$, $\xi^+(t)=\max\{t^{p^-},t^{p^+}\}$, for all $t\geq0$.
    Then
    \begin{itemize}
      \item[(i)] $\xi^-(\|u\|_{L^G(\R^n)})\leq\Phi_{G,\R^n}(u)\leq\xi^+(\|u\|_{L^G(\R^n)})\ \text{for}\ u\in L^{G}(\R^n)$,
      \item[(ii)] $\xi^-([u]_{W^{s,G}(\R^n)})\leq \Phi_{s,G,\R^n}(u) \leq\xi^+([u]_{W^{s,G}(\R^n)})\ \text{for}\ u\in W^{s,G}(\R^n)$.
    \end{itemize}
 \end{lema}

\section{Lagrange multipliers and the eigenvalue problem} \label{section.lagrange}
In this section we define the notion of Dirichlet, Neumann and Robin eigenvalue problems in the context of fractional Orlicz-Sobolev spaces. We recall some existence results already proved in \cite{VS20} and \cite{S19} for the Dirichlet case, and state further extension to more general boundary conditions.

 We say that
\begin{itemize}
\item[] $\lambda$ is an \emph{eigenvalue} of \ref{eq.d} with \emph{eigenfunction} $u\in W_0^{s,G}(\Omega)$ if
\begin{equation}\label{weak.eq.d}
  \langle (-\Delta_g)^s u,v \rangle=\lambda\int_{\Omega} g(|u|)\frac{u}{|u|}v\, dx\qquad \forall v\in W_0^{s,G}(\Omega).
\end{equation}

\item[] $\lambda$ is an \emph{eigenvalue} of \ref{eq.n} with \emph{eigenfunction} $u\in W^{s,G}_*(\Omega)$ if
\begin{equation*} \label{eig.n}
\langle (-\Delta_g)^s u,v\rangle_*  = \lambda\int_\Omega g(|u|)  \frac{u}{|u|} v \,dx \qquad \forall v \in W^{s,G}_*(\Omega).
\end{equation*}

\item[] $\lambda$ is an \emph{eigenvalue} of \ref{eq.r} with  with \emph{eigenfunction} $u\in \mathcal{X}_\beta(\Omega)$  if
\begin{equation*} \label{eig.r}
\langle (-\Delta_g)^s u,v \rangle_* = \lambda\int_{\Omega} g(|u|)\frac{u}{|u|}v\, dx - \int_{\mathbb{R}^n\setminus \Omega} \beta g(|u|)\frac{u}{|u|}v\, dx \qquad \forall v\in \mathcal{X}_\beta(\Omega).
\end{equation*}
\end{itemize}


As mentioned, eigenvalues for non-homogeneous eigenproblems strongly depend on the energy level $\alpha>0$. It is possible to prove existence of an eigenvalue by a multiplier argument once it is proved existence   of the following associated minimization problems
\begin{align*}
\Lambda^D_\alpha &:=\min\{ \mathcal{J}(u) \colon u\in M^D_\alpha\}, \quad \text{ where } M^D_\alpha=\{ u \in W^{s,G}_0(\Omega) \colon \mathcal{I}(u)=\alpha\}\\
\Lambda^N_\alpha &=\min\{ \mathcal{J}_0(u) \colon u\in M^N_\alpha\},\quad \text{ where } M^N_\alpha =\{ u \in W^{s,G}_*(\Omega) \colon \mathcal{I}(u)=\alpha\}\\
\Lambda^R_\alpha &:=\min\{ \mathcal{J}_\beta(u) \colon u\in M^R_\alpha\}, \quad \text{ where } M^R_\alpha =\{ u \in \mathcal{X}_\beta(\Omega) \colon \mathcal{I}(u)=\alpha\}.
\end{align*}

As pointed out in \cite{BS20}, by standard computations of the calculus of variations, for each $\alpha>0$ the quantities $\Lambda_\alpha^D$, $\Lambda_\alpha^N$ and $\Lambda_\alpha^R$ are attained by  suitable functions $u_\alpha^D$, $u_\alpha^N$ and $u_\alpha^R$, respectively. Since condition \eqref{G1} is assumed on $G$, the functionals $\mathcal{J}$, $\mathcal{J_\beta}$  and $\mathcal{I}$ are class $C^1$ with Frech\'et derivatives given by
\begin{align} \label{derivadas}
\begin{split}
\langle (\mathcal{J}'(u),v \rangle &= \langle (-\Delta_g)^su,v\rangle \qquad  \forall v \in W^{s,G}_0(\Omega),\\
\langle (\mathcal{J}_0'(u),v \rangle &= \langle (-\Delta_g)^su,v\rangle_* \qquad \forall v \in W^{s,G}_*(\Omega),\\
\langle (\mathcal{I}'(u),v \rangle &= \int_\Omega g(|u|)\frac{u}{|u|}v\,dx \qquad \forall v \in L^G(\Omega),\\
\langle (\mathcal{J}'_\beta(u),v \rangle &= \langle (-\Delta_g)^su,v\rangle_* + \int_{\R^n\setminus \Omega} \beta g(|u|)\frac{u}{|u|}v\,dx \qquad \forall v \in W^{s,G}_*(\Omega),
\end{split}
\end{align}
and therefore, by an application of the Lagrange multipliers rule, there exist numbers $\lambda^D_\alpha$, $\lambda^N_\alpha$ and $\lambda^R_\alpha$ which are eigenvalues of \ref{eq.d}, \ref{eq.n} and \ref{eq.r} with eigenfunctions $u_\alpha^D$, $u_\alpha^N$ and $u_\alpha^R$, respectively. Furthermore, as proved in \cite[Proposition 5.6 and 5.8]{BS20}, minimizers  are comparable each other; more precisely, $\Lambda_\alpha^N\leq \Lambda_\alpha^R \leq \Lambda_\alpha^D$, and a similar relation for the eigenvalues holds, up to a multiplicative constant.

Proceeding as in \cite[Theorem 1.3]{VS20} (see Theorem \ref{lam.inf}) it can be proved that for each fixed $\alpha_0>0$ it holds that
\begin{align*}
\hat \lambda^D &:=\inf\{ \lambda^D_\alpha\colon \alpha \in (0,\alpha_0)\},\\
\hat \lambda^N &:=\inf\{ \lambda^N_\alpha\colon \alpha \in (0,\alpha_0) \},\\
\hat \lambda^R &:=\inf\{ \lambda^R_\alpha\colon \alpha \in (0,\alpha_0) \}
\end{align*}
are eigenvalues of \ref{eq.d}, \ref{eq.n} and \ref{eq.r}, respectively.

In contrast with $p-$Laplacian type problems, when dealing with non-homogeneous eigenproblems, eigenvalues are not variational in general. One could consider the variational quantity
\begin{equation*}
\bar \lambda^D=\displaystyle\inf_{u\in W^{s,G}_0(\Omega)\setminus\{0\}}\frac{\iint_{\mathbb{R}^n \times\R^n } g(|D_s u|)|D_s u|\,d\mu}{\int_{\Omega} g(|u|)|u|\, dx},
\end{equation*}
but in general this number cannot be probed to be an eigenvalue. In an analogous way, we can define $\bar \lambda^R$ and $\bar \lambda^N$, and the same assertion holds. However, the following result is true:

\begin{thm}
It holds that $\bar \lambda^D \leq \hat \lambda^D$ and there is no eigenvalue $\lambda$ of \ref{eq.d} such that $\lambda< \bar \lambda^D$.

The same holds by changing the superscript $D$ with $N$ or $R$.
\end{thm}

\begin{proof}
First observe that, since $\hat \lambda^D$ is an eigenvalue of \ref{eq.d}, there exists a nontrivial eigenfunction $\hat u \in W^{s,G}_0(\Omega)$ such that \eqref{weak.eq.d} holds, therefore
$$
\lambda_0 =\inf\left\{ \frac{\langle (-\Delta_g)^s   u,  u\rangle}{\int_\Omega g(|u|)|u|\,dx}\colon u\in W^{s,G}_0(\Omega)\right\} \leq \frac{\langle (-\Delta_g)^s \hat u,\hat u\rangle}{\int_\Omega g(|u|)|u|\,dx} = \hat \lambda^D.
$$

If we suppose that there exists  $\lambda< \lambda_0$ which is an eigenvalue of problem \ref{eq.d} with eigenfunction $u_\lambda\in W^{s,G}_0(\Omega)$, we arrive to a contradiction since
$$
\lambda =  \frac{\langle (-\Delta_g)^s   u_\lambda,\hat u_\lambda\rangle}{\int_\Omega g(|u_\lambda|)|u_\lambda|\,dx} < \lambda_0 = \inf\left\{ \frac{\langle (-\Delta_g)^s   u,  u\rangle}{\int_\Omega g(|u|)|u|\,dx}\colon u\in W^{s,G}_0(\Omega)\right\}.
$$

The proofs of the Robin and Neumann case run analogously.
\end{proof}

We finish this section by proving Theorem \ref{lam.inf}.

\begin{proof}[Proof of Theorem \ref{lam.inf}]
Let  $\{\lambda_{\alpha_k}\}_{k\in\N}$ be a sequence of eigenvalues of \ref{eq.d} such that $\lambda_{\alpha_k}\to \lambda$ and let $\{u_{\alpha_k}\}_{k\in\N}\subset W^{s,G}_0(\Omega)$ be the corresponding sequence of associated eigenfunctions such that $\Phi_{G,\Omega}(u_{\alpha_k})=\alpha_k$ with $\alpha_k \in (0,\alpha_0)$ for all $k\in\N$. Then, it holds that
\begin{equation} \label{eqj}
\iint_{\R^n\times\R^n} g(|D_s u_{\alpha_k}|) \frac{D_s u_{\alpha_k}}{|D_s u_{\alpha_k}|}  D_s v\,d\mu = \lambda_{\alpha_k} \int_\Omega g(|u_{\alpha_k}|)\frac{u_{\alpha_k}}{|u_{\alpha_k}|} v \qquad \forall v\in W^{s,G}_0(\Omega).
\end{equation}

Observe that $\{u_{\alpha_k}\}_{k\in\N}$ is bounded in $W^{s,G}_0(\Omega)$ since by \eqref{G1} and \eqref{eqj}
\begin{align*}
\Phi_{s,G,\R^n}(u_{\alpha_k}) &\leq \frac{1}{p^-}\iint_{\R^n\times\R^n} g(|D_s u_{\alpha_k}|) |D_s u_{\alpha_k}| \,d\mu\\
&= \frac{\lambda_{\alpha_k}}{p^-} \int_\Omega g(|u_{\alpha_k}|)|u_{\alpha_k}|  \leq \frac{p^+}{p^-} \lambda_{\alpha_k} \alpha_k\\
&\leq \frac{p^+}{p^-} (1+\lambda) \alpha_0
\end{align*}
for $k$ big enough. Then, by using the compact embedding given in Proposition \ref{ceb} (i), up to a subsequence,  there exists $u\in W^{s,G}_0(\Omega)$ such that
\begin{equation}\label{convergen}
\begin{array}{ll}
u_{\alpha_k}\to u &\text{ strongly  in }L^{G}( \Omega),\\
u_{\alpha_k}\to  u &\text{ a.e. in } \R^n.
\end{array}
\end{equation}
From the continuity of $t\mapsto g(t)\frac{t}{|t|}$ and \eqref{convergen} we deduce that
$$
g(|D_s u_k|)\frac{D_s u_k}{|D_s u_k|}\to g(|D_s u|)\frac{D_s u}{|D_s u|}\quad  \text{ a.e. in } \Omega
$$
and also from \eqref{convergen}
$$
\limsup_{k\to\infty} \Phi_{G,\R^n}(u_{\alpha_k}) \in (0,\alpha_0).
$$
Hence, taking limit as $k\to\infty$ in \eqref{eqj} we obtain that
$$
\iint_{\R^n\times\R^n} g(|D_s u|)\frac{D_s u}{|D_s u|}   D_s v \, d\mu = \lambda \int_\Omega g(|u|)\frac{u}{|u|} v \qquad \text{ for all } v\in W^{s,G}_0(\Omega)
$$
from where the proof in the Dirichlet case follows.

The proof for the Neumann and Robin case follows analogously by replacing $W^{s,G}_0(\Omega)$ with $W^{s,G}_*(\Omega)$ or $\mathcal{X}_\beta(\Omega)$, respectively, and using the compact embedding given in Proposition \ref{ceb} item (ii) to deduce \eqref{convergen}.
\end{proof}

\section{Ljusternik-Schnirelman eigenvalues} \label{LS}

In order to build a sequence of eigenvalues for problems \ref{eq.d}, \ref{eq.n} and \ref{eq.r}, the idea is to apply an abstract theorem of the so-called Ljusternik-Schnirelman theory. See for instance \cite{Br65, Papageorgiou, zeidlerIII}. We  will particularly use  the result stated in \cite[Theorem 9.27]{Papageorgiou}.

Given $\alpha>0$, assume that $\mathcal{A}$, $\mathcal{B}$ are two functionals defined in a reflexive Banach space $\mathcal{X}$, such that
\begin{itemize}
\item[($h_1$)] $\mathcal{A}$,$\mathcal{B}$  are  $C^1(\mathcal{X},\R)$ even functionals  with $\mathcal{I}(0)=\mathcal{B}(0)=0$ and the level set
$$
M_\alpha := \{ u\in \mathcal{X} \colon \mathcal{B}(u)=\alpha\}
$$
is bounded.

\item [($h_2$)] $\mathcal{A}'$ is strongly continuous, i.e.,
$$
u_k\rightharpoonup u \text{ in } \mathcal{X} \implies \mathcal{A}'(u_k) \to \mathcal{A}'(u).
$$
Moreover, for any $u$ in the closure of the convex hull of $M_\alpha$,
$$
\langle \mathcal{A}'(u),u\rangle=0 \iff \mathcal{A}(u)=0 \iff u=0.
$$

\item [($h_3$)] $\mathcal{B}'$ is continuous, bounded and, as $k\to\infty$, it holds that
$$
u_k\rightharpoonup u, \quad \mathcal{B}'(u_k)\rightharpoonup v, \quad \langle \mathcal{B}'(u_k),u_k\rangle\to\langle v,u\rangle \implies u_k\to u \text{ in } \mathcal{X}.
$$

\item [($h_4$)] For every $u \in \mathcal{X} \setminus\{0\}$ it holds that
$$
\langle \mathcal{B}'(u),u\rangle>0,\qquad \displaystyle\lim_{t\to+\infty}\mathcal{B}(tu)=+\infty,\qquad \displaystyle\inf_{u\in M_\alpha}\langle \mathcal{B}'(u),u\rangle>0.
  $$
\end{itemize}

Define max-min values
\begin{equation*}\label{ak}
  c_{k,\alpha}=\begin{cases}
        \sup_{K\in \mathcal{C}_k}\inf_{u\in K}\mathcal{A}(u), & \mathcal{C}_k\neq\emptyset , \\
        0, &  \mathcal{C}_k= \emptyset,
      \end{cases}
\end{equation*}
where, for any $k\in\mathbb{N}$,
$$
\mathcal{C}_k:=\{ K\subset M_\alpha \text{ compact, symmetric with } \mathcal{A}(u)>0 \text{ on } K  \text{ and } \gamma(K)\geq k\},
$$
and the Krasnoselskii genus of $K$ is defined as
$$
\gamma(K):=\inf\{p\in\mathbb{N}\colon \exists h\colon K\to \R^p\setminus\{0\}\ \text{such that}\ h \text{ is continuous and odd} \}
$$
see \cite{K64} for details.

Thus,  $\{c_{k,\alpha}\}_{k\geq1}$ forms a nonincreasing sequence
$$
+\infty> c_{1,\alpha}\geq c_{2,\alpha}\geq \ldots\geq c_{k,\alpha}\geq\ldots\geq0.
$$

Under these considerations, the Ljusternik–Schnirelmann principle stated in \cite[Theorem 9.27]{Papageorgiou} establishes that there exists a sequence $\{(\mu_{k,\alpha},u_{k,\alpha})\}_{k\geq1}$ such that
$$
  \langle \mathcal{A}'(u_{k,\alpha}),v\rangle=\mu_{k,\alpha} \langle \mathcal{B}'(u),v\rangle \qquad \forall v\in \mathcal{X}^*
$$
such that $u_{k,\alpha} \in M_\alpha$, $\mathcal{A}(u_{k,\alpha})=c_{k,\alpha}$,  $\mu_{k,\alpha}\neq 0$, $\mu_{k,\alpha}\to 0$, and $u_{k,\alpha}\rightharpoonup 0$ in $\mathcal{X}$.

\subsection{The Dirichlet case}

With these preliminaries, we are in position to prove  Theorem \ref{main.result.1}.

We consider the  space $\mathcal{X}:=W^{s,G}_0(\Omega)$ and the functionals $\mathcal{B}(u):=\mathcal{J}(u)$ and $\mathcal{A}(u):=\mathcal{I}(u)$ defined in \eqref{funcionales}. As mentioned, these functionals are $C^1$ and their Frech\'et derivatives are given in \eqref{derivadas}.

\begin{lema} \label{hipotesis}
The functionals $\mathcal{J}$ and $\mathcal{I}$ defined above   fulfill hypotheses $(h_1)$--$(h_4)$.
\end{lema}
\begin{proof}

$(i)$ Clearly, the maps $\mathcal{I},\mathcal{J}$ are even and $\mathcal{I}(0)=\mathcal{J}(0)=0$.

 \medskip

 $(ii)$ We notice that
 $$
p^- \xi^-(\|u\|_{L^G(\Omega)})\leq p^- \mathcal{I}(u)\leq \langle\mathcal{I}'(u),u\rangle \leq p^+ \mathcal{I}(u)\leq p^+\xi^+(\|u\|_{L^G(\Omega)}).
$$
Then immediately we obtain
$$
\langle \mathcal{I}'(u),u\rangle=0 \ \Leftrightarrow \ \mathcal{I}(u)=0 \ \Leftrightarrow \ u=0.
$$
 Thus, it remains to check that $\mathcal{I}'$ is strongly continuous. Let $u_k\rightharpoonup u$ in $W^{s,G}_0(\Omega)$, then $\{u_k\}_{k\in\mathbb{N}}$ is bounded in $W^{s,G}_0(\Omega)$. We need to show that $\mathcal{I}'(u_k)\to \mathcal{I}'(u)$ in $W^{-s,\tilde G}(\Omega)$.
  \begin{align*}
	|\langle \mathcal{I}'(u_{k})-\mathcal{I}'(u),v\rangle|&=\left|\int_{\Omega}\left(g(|u_k|)\frac{u_k}{|u_{k}|}-g(|u|)\frac{u}{|u|}\right) v\,dx\right|\\
	&\leq\left|\int_{\Omega} g(|u_k|)\left(\frac{u_{k}}{|u_{k}|}-\frac{u}{|u|}\right) v\,dx\right| +\left|\int_{\Omega}(g(|u_k|)-g(|u|))\frac{u}{|u|}v\,dx\right| \\
	&:=I_{1,k}+I_{2,k}.
	\end{align*}
By using  Lemma \ref{lemita.1} and H\"{o}lder's inequality, the first term can be bounded as
$$
I_{1,k}\leq \|g(|u_k|)\|_{\tilde G}\|2 v\|_{G}\rightarrow 0,\ k\rightarrow+\infty.
$$
Similarly, by using the $\Delta_2$ condition, Lemma \ref{lemita.1} and H\"{o}lder's inequality we get
$$
I_{2,k}\leq \|g(|u_k|)-g(|u|)\|_{\tilde G}\|v\|_{G}.
$$
Since $u_k \rightharpoonup u$ in $W^{s,G}_0(\Omega)$, in light of Proposition \ref{ceb}, $u_k\to u$ strongly in $L^G(\Omega)$ and a.e. in $\R^n$, moreover, from Lemma \ref{lemita}, $g(|u_k|)-g(|u|) \in L^1(\Omega)$, therefore, by dominated convergence theorem, $\int_\Omega \tilde G(|g(|u_k|)-g(|u|)|)\,dx \to 0$ and hence $\|g(|u_k|)-g(|u|)\|_{L^{\tilde G}(\Omega)} \to 0$ as $k\to\infty$.

From the last relations it follows that $\|\mathcal{I}'(u_{k})-\mathcal{I}'(u)\|_{W^{-s,\tilde G}(\Omega)}\rightarrow0$ as required.

\medskip

$(iii)$ One can easily see that $\mathcal{J}'$ is continuous (see for instance \cite[Proposition 4.1]{S19}). From Lemma \ref{lemita.1} and Holder's inequality
$$
|\langle \mathcal{J}'(u),v| \rangle \leq \|g(|u|)\|_{L^{\tilde G}(\R^n\times\R^n,d\mu)} \|v\|_{L^{G}(\R^n\times\R^n,d\mu)}\leq p^+[u]_{W^{s,G}(\R^n)} [v]_{W^{s,G}(\R^n)}
$$
from there $\mathcal{J}'$ is bounded.

It remains to be showed that if $\{u_k\}_{k\in\mathbb{N}}$ is a sequence in $W^{s,G}_0(\Omega)$ such that
\begin{equation} \label{asump}
  u_k\rightharpoonup u, \quad \mathcal{J}'(u_k)\rightharpoonup v, \quad \langle \mathcal{J}'(u_k),u_k\rangle\to\langle v,u\rangle\
\end{equation}
then $u_k \to u$ in $W^{s,G}_0(\Omega)$.

Since $G$ is convex, we have
$$
G(|D_s u|)\leq G\bigg{(}\bigg{|}\frac{D_s u+D_s u_k}{2}\bigg{|}\bigg{)}+g(|D_s u|)\frac{D_s u}{|D_s u|}.\frac{D_s u-D_s u_k}{2}
$$
and
$$
G(|D_s u_k|)\leq G\bigg{(}\bigg{|}\frac{D_s u+D_s u_k}{2}\bigg{|}\bigg{)}+g(|D_s u_k|)\frac{D_s u_k}{|D_s u_k|}.\frac{D_s u_k-D_s u}{2}
$$
for every $u,v\in W^{s,G}_0(\Omega)$. Adding the above two relations and integrating over $\R^n$ we find that
\begin{align}\label{umo1}
 \frac{1}{2}\iint_{\R^n\times \R^n} \Big(g(|D_s u|)\frac{D_s u}{|D_s u|}&-g(|D_s u_k|)\frac{D_s u_k}{|D_s u_k|}  \Big)(D_s u-D_s u_k)\ d\mu\nonumber\\
 &\geq \Phi_{s,G,\R^n}(u) + \Phi_{s,G,\R^n}(u_k) - 2 \Phi_{s,G,\R^n}\left( \frac{u+u_k}{2}\right)
\end{align}
for every $u,v\in W^{s,G}_0(\Omega)$. By applying Lemma \ref{lemita}, the right part of the inequality above can be bounded by below by $2\Phi_{s,G,\R^n}\left(\frac{u-u_k}{2}\right)$, and hence we get
$$
 \langle (-\Delta_g)(u)- (-\Delta_g)(u_k), u-u_k \rangle \geq 4 \Phi_{s,G,\R^n}(u-u_k).
$$
This, together with Lemma \ref{ineq1} yields
\begin{equation} \label{eqxx1}
\langle \mathcal{J}'(u)-\mathcal{J}'(u_k), u-u_k\rangle \geq \xi^-([u-u_k]_{W^{s,G}(\R^n)}).
\end{equation}

On the other hand, Proposition \ref{ceb} gives that $u_k\to u$ in $L^G(\Omega)$ and a.e. in $\R^n$, which, mixed up with the assumptions \eqref{asump} allows us to deduce that
$$
\displaystyle\lim_{k\to\infty}\langle \mathcal{J}'(u_k)-\mathcal{J}'(u),u_k-u \rangle=\displaystyle\lim_{k\to\infty} \left(\langle \mathcal{J}'(u_k),u_k\rangle-\langle \mathcal{J}'(u_k),u\rangle-\langle \mathcal{J}'(u),u_k-u\rangle \right)=0.
$$
Hence, from \eqref{eqxx1}, $[u_k-u]_{W^{s,G}(\R^n)}\to0$ as $k\to\infty$ as required.

\medskip

$(iv)$ It is clear that, for any $u\in W^{s,G}_0(\Omega)\setminus\{0\}$,
$$
\langle \mathcal{J}'(u),u\rangle>0,\qquad \displaystyle\lim_{t\to+\infty}\mathcal{J}(tu)=+\infty,\qquad \displaystyle\inf_{u\in M_\alpha}\langle \mathcal{J}'(u),u\rangle>0.
  $$
This concludes the proof.
\end{proof}

\begin{proof}[Proof of Theorem \ref{main.result.1} (Dirichlet case)]
In light of Lemma \ref{hipotesis} and in virtue of  \cite[Theorem 9.27]{Papageorgiou}, there exist a sequence of positive numbers $\{\mu_k\}_{k\in\mathbb{N}}$ tending to 0 and a corresponding sequence of functions $\{u_{k,\alpha}^D\}_{k\in\N}\in W^{s,G}_0(\Omega)$ such that
$$
\langle (-\Delta_g)^s u_{k,\alpha}^D,v \rangle=\frac{1}{\mu_{k,\alpha}^D}\int_{\Omega} g(|u_{k,\alpha}^D|)\frac{u_{k,\alpha}^D}{|u_{k,\alpha}^D|}v\, dx \quad \forall v\in W^{s,G}_0(\Omega).
$$
Moreover, $\mathcal{J}(u_{k,\alpha}^D)=\alpha$ and
$$
\mathcal{I}(u_{k,\alpha}^D):=c_{k,\alpha}^D=\sup_{K\in \mathcal{C}_k}\inf_{u\in K}\mathcal{I}(u)>0.
$$
Consequently, $\lambda_{k,\alpha}^D=1/\mu_{k,\alpha}$ is an eigenvalue of \ref{eq.d} with eigenfunction $u_{k,\alpha}^D$.
\end{proof}

\begin{rem}
We mention some observations regarding the Ljusternik-Schnirelman sequence of eigenvalues obtained in Theorem \ref{main.result.1}.
\begin{itemize}
\item[(a)] Since in general these functionals are not homogeneous, we cannot claim that $\lambda_{k,\alpha}^D=1/c_{k,\alpha}^D$.

\item[(b)] The first Ljusternik-Schnirelman eigenvalue $\lambda_{1,\alpha}^D$ is in fact different to the eigenvalue $\lambda^D_\alpha$ obtained in Section \ref{section.lagrange} by means of the Lagrange multipliers rule.

\item[(c)] The sequence $\{\lambda^D_{k,\alpha}\}_{k\in\mathbb{N}}$  does not exhaust the spectrum of \ref{eq.d}. In fact, the spectrum is not discrete since  the parameter $\alpha$ can be taken in $\R^+$.

\item[(d)] The eigenfunction associated to $\lambda_{k,\alpha}^D$ is one-signed since the $\Phi_{s,G,\R^n}(\cdot)$ and $\Phi_{G,\Omega}(\cdot)$ are invariant by replacing $u_{k,\alpha}^D$ with $|u_{k,\alpha}^D|$.

\item[(e)] As mentioned in Section \ref{section.lagrange}, due to the closedness of the spectrum $\Sigma_D$ we have that the quantity $\inf\{ \lambda_{k,\alpha}^D\colon 0<\alpha<\alpha_0\}$ is also an eigenvalue of \ref{eq.d}, for any $\alpha_0>0$ fixed.
\end{itemize}

\end{rem}

\subsection{The Neumann/Robin case}

To deal with the Robin case we take $\mathcal{X}:=\mathcal{X}_\beta(\Omega)$, $\mathcal{A}:=\mathcal{I}(u)$ and $\mathcal{B}:=\mathcal{J}_\beta(u)$. The Neumann case it follows just by setting $\beta=0$.

\begin{proof}[Proof of Theorem \ref{main.result.1} (Neumann/Robin case)]
  The proof of this result is similar to the proof of the Dirichlet case, just noticing that, the embedding $\mathcal{X}_\beta(\Omega)\hookrightarrow L^{G}(\Omega)$ is compact (see Proposition \ref{ceb} (ii)) and the quantities $[u]_{W^{s,G}(\R^n)}$ and $[u]_{W^{s,G}_*(\R^n)}$ play a symmetrical role.
\end{proof}

\section{Minimax eigenvalues}

This section is devoted to prove Theorem \ref{main.result.2}. We start with the Dirichlet case.

Given $\alpha>0$ we recall that
$$
M_\alpha^D:=\{u\in W^{s,G}_0(\Omega)\colon \mathcal{I}(u)=\alpha\}
$$
defines a $C^1$ manifold. We denote $\tau_\alpha$ the restriction of $(-\Delta_g)^s$ to $M_\alpha^D$, i.e., for each $u\in M_\alpha^D$,
$$
\langle \tau_\alpha(u), v\rangle :=\langle (-\Delta_g)^s u,v\rangle \quad \forall v\in T_u M_\alpha^D
$$
where the tangent space of $M_\alpha^D$ at $u$ is defined as
$$
T_u M_\alpha^D :=\left\{ v\in W^{s,G}_0(\Omega) \colon  \langle \mathcal{I}'(u),v\rangle =0 \right\}.
$$
Thus, we define
$$
\|\tau_\alpha u\|:= \| \tau_\alpha u\|_{(T_u M_\alpha^D)^*}
$$
where $(T_u M_\alpha^D)^*$ is the dual space of $T_u M_\alpha^D$.

We recall that the \emph{duality mapping} $J\colon W^{s,G}_0(\Omega) \to (W^{s,G}_0(\Omega))^*$ (where we have denoted the dual space of the space $\mathcal{X}$ as $\mathcal{X}^*$) is defined as a bijective isometry such that
\begin{equation} \label{relacJ}
\|Ju\|_{(W^{s,G}_0(\Omega))^*} = [u]_{W^{s,G}_0(\Omega)}, \qquad \langle Ju, u \rangle = [u]_{W^{s,G}(\R^n)}^2 \qquad \forall u\in W^{s,G}_0(\Omega).
\end{equation}
See \cite[Chapter 18.11c]{zeidlerIIA} and \cite[Proposition 47.18]{zeidlerIII} for details.

Given $\alpha>0$, the idea  is to apply the minimax theorem stated in \cite{Cuesta} to the functional $\mathcal{J}(u)$ under the constraint $\mathcal{I}(u)=\alpha$ to obtain a sequence of critical points of the form
$$
C_{k,\alpha}^D = \inf_{h\in \Gamma(S^{k-1},M_\alpha^D)} \sup_{w \in S^{k-1}} \Phi_{s,G,\R^n}(h(w)).
$$
being $\Gamma(S^{k-1}k,M_\alpha^D)=\{h \in C(S^{k-1},M_\alpha^D)\colon h \text{ is odd}\}$.

Recall that the derivatives of $\mathcal{I}$ and $\mathcal{J}$ are given by  $\mathcal{J}'(u)=(-\Delta_g)^s u$ and $\mathcal{I}'(u)=\frac{g(|u|) u}{|u|}$ for any $u\in L^G(\Omega)$.

\begin{defn}
We say that $\mathcal{J}$ satisfies the \emph{Palais-Smale condition on $M_\alpha^D$ at level $c$} if any sequence $\{u_n\}_{n\in\N}\subset M_\alpha^D$ such that $\mathcal{J}(u_n)\to c$ and $\|\tau_\alpha u_n\| \to 0$, possesses a convergent subsequence.
\end{defn}

The key ingredient is to analyze the validity of the Palais-Smale condition.

\begin{lema} \label{lemaPS}
Given $\alpha>0$ the functional $\mathcal{J}$ satisfies the Palais-Smale condition on $M_\alpha^D$ at level $C_{k,\alpha}^D$.
\end{lema}
\begin{proof}
We follow closely the construction of \cite[Theorem 5.3]{szulkin}. Given $k\in\N$, let $\{u_n\}_{n\in\N}$ be a sequence in $M_\alpha^D$ such that
\begin{equation} \label{hyp}
\mathcal{I}(u_n)\to C_{k,\alpha}^D \quad \text{ and } \quad \|\tau_\alpha (u_n)\| \to 0.
\end{equation}
For each $u\in M_\alpha^D$ define the projection $P_u\colon W^{s,G}_0(\Omega) \to T_u M_\alpha^D$ such that
$$
P_u v = v-\frac{\langle \mathcal{I}'(u),v \rangle } {\|\mathcal{I}'(u)\|^2_{(W^{s,G}_0(\Omega))^*}} J^{-1}\left(  \mathcal{I}'(u) \right).
$$
Observe first that $P_u$ is well defined: given $v\in W^{s,G}_0(\Omega)$ and $u\in M_\alpha^D$ we have
$$
\langle \mathcal{I}'(u),P_u v \rangle = \langle \mathcal{I}'(u),v \rangle -\frac{\langle \mathcal{I}'(u),v \rangle } {\|\mathcal{I}'(u)\|^2_{(W^{s,G}_0(\Omega))^*}} \langle \mathcal{I}'(u),  J^{-1}\left(  \mathcal{I}'(u) \right) \rangle=0
$$
since due to \eqref{hyp}, $\langle \mathcal{I}'(u),  J^{-1}\left(  \mathcal{I}'(u) \right) \rangle= \| \mathcal{I}'(u)\|^2_{(W^{s,G}_0(\Omega))^*}$.

Moreover,  for every $v\in W^{s,G}_0(\Omega)$ we get the inequality
\begin{align*}
\langle \mathcal{I}'(u),v \rangle  \leq \|\mathcal{I}'(u)\|_{(W^{s,G}_0(\Omega))^*} [v]_{W^{s,G}(\R^n)}.
\end{align*}
and by \eqref{relacJ} it follows that
$$
\left[J^{-1}\left(  \mathcal{I}'(u) \right) \right]_{W^{s,G}(\R^n)} = \|\mathcal{I}'(u)\|_{ (W^{s,G}_0(\Omega))^* }.
$$
Therefore we get that $[P_u v]_{W^{s,G}(\R^n)} \leq 2 [v]_{_{W^{s,G}(\R^n)}}$, which implies that, for any $v\in W^{s,G}_0(\Omega)$
$$
\langle (-\Delta_g)^s u, P_u v\rangle = \langle \tau_\alpha(u), P_u v\rangle  \leq 2\| \tau_\alpha(u)\| [v]_{W^{s,G}(\R^n)}.
$$
Consequently, by \eqref{hyp}
$$
\sup_{[v]_{W^{s,G}(\R^n)}\leq 1} |\langle (-\Delta_g)^s u_n, P_{u_n} v\rangle| \leq  \sup_{[v]_{W^{s,G}(\R^n)}\leq 1} 2 \| \tau_\alpha(u_n)\| [v]_{W^{s,G}(\R^n)}\to 0 \quad \text{as } n\to\infty.
$$
Hence we get that
\begin{align*}
\langle (-\Delta_g)^s  u_n, P_{u_n} v \rangle &= \langle (-\Delta_g)^s u_n,v \rangle- \left\langle (-\Delta_g)^s u_n, \frac{ \langle \mathcal{I}' u_n ,v \rangle }{\|\mathcal{I}'(u)\|^2_{(W^{s,G}_0(\Omega))^*}}  J^{-1}(\mathcal{I}'(u_n)) \right\rangle\\
&= \langle (-\Delta_g)^s u_n,v \rangle- \left\langle (-\Delta_g)^s u_n,   J^{-1}(\mathcal{I}'(u_n)) \right\rangle \frac{ \langle \mathcal{I}'(u_n) ,v \rangle }{ \|\mathcal{I}'(u)\|^2_{(W^{s,G}_0(\Omega))^*} } \\
&= \langle (-\Delta_g)^s u_n,v \rangle- \left\langle \frac{\left\langle (-\Delta_g)^s u_n,   J^{-1}(\mathcal{I}'(u_n)) \right\rangle }{ \|\mathcal{I}'(u)\|^2_{(W^{s,G}_0(\Omega))^*}} \mathcal{I}'(u_n) ,v \right\rangle \to 0,
\end{align*}
that is
$$
(-\Delta_g)^s u_n -  \frac{\left\langle (-\Delta_g)^s u_n,   J^{-1}(\mathcal{I}'(u_n)) \right\rangle }{  \|\mathcal{I}'(u)\|^2_{(W^{s,G}_0(\Omega))^*} } \mathcal{I}'(u_n) \rightharpoonup 0 \quad \text{ weakly in }W^{s,G}_0(\Omega).
$$

From \eqref{hyp}, up to a subsequence $u_n \rightharpoonup u$ weakly in $W^{s,G}_0(\Omega)$ and strongly in $L^G(\Omega)$ to some $u\in W^{s,G}_0(\Omega)$ due to Proposition \ref{ceb}. Observe that in Lemma \ref{hipotesis} we have proved that $\mathcal{I}'$ satisfies property $(h_2)$ of Section \ref{LS}, i.e., $\mathcal{I}$ is strongly continuous, which implies that $\mathcal{I}'(u_n)\to \mathcal{I}'(u)$.

Moreover, \eqref{hyp} also gives that $(-\Delta_g)^s u_n$ is bounded and $u_n$ is bounded away from zero. Therefore, up to a subsequence
\begin{equation} \label{desJ}
(-\Delta_g)^s u_n \rightharpoonup v \quad \text{ weakly in } (W^{s,G}_0(\Omega))^*
\end{equation}

In Lemma \ref{hipotesis} we have proved that $\mathcal{J}' = (-\Delta_g)^s$ satisfies property $(h_3)$ of Section \ref{LS}, from where, in view of \eqref{desJ} we have that $u_n \to u$ strongly in $W^{s,G}_0(\Omega)$. This proves that $\mathcal{J}$ satisfies the Palais-Smale condition on $M_\alpha^D$ at level $C_{k,\alpha}^D$ which concludes the proof.
\end{proof}

\begin{rem}
Observe that in the proof of Lemma \ref{lemaPS} we have not used any particular property for functions in $W^{s,G}_0(\Omega)$. With the pertinent changes, the same arguments can be applied to deduce that, given $\alpha>0$, the functional $\mathcal{J}$ satisfies the Palais-Smale condition on $M_\alpha^N$ (resp. $M_\alpha^R$) at level $C_{k,\alpha}^N$ (resp. $C_{k,\alpha}^R$), where $C_{k,\alpha}^N$ (resp. $C_{k,\alpha}^R$) is defined just by changing the superscript $D$ by $N$ (resp. $R$) in the definition of $C_{k,\alpha}^D$. We leave to the reader the remaining details.
\end{rem}


\begin{proof}[Proof of Theorem \ref{main.result.2} ]

Due to Lemma \ref{lemaPS}, $\mathcal{J}$ (resp. $\mathcal{J}_0$ ,$\mathcal{J}_{\beta}$) satisfies the Palais-Smale condition on $M_\alpha^D$
(resp. $M_\alpha^N$, $M_\alpha^R$) at level $C_{\alpha,k}^D$ (resp. $C_{\alpha,k}^N$, $C_{\alpha,k}^R$) for each $k\in\mathbb{N}$, then by
\cite[Proposition 2.7]{Cuesta}
 there exists $u_{\alpha,k}^D\in M_\alpha^D$ (resp. $u_{\alpha,k}^N\in M_\alpha^N$, $ u_{\alpha,k}^R\in M_\alpha^R$) such that
$$
\mathcal{J}(u_{k,\alpha}^D)=C_{k,\alpha}^D, \qquad (\text{resp.}\ \mathcal{J}_0(u_{k,\alpha}^N)=C_{k,\alpha}^N,\,\,  \mathcal{J}_\beta (u_{k,\alpha}^R)=C_{k,\alpha}^R)
$$
and
$$
\mathcal{J}'(u_{k,\alpha}^D)=\mathcal{J}_0'(u_{k,\alpha}^N)=\mathcal{J}_\beta'(u_{k,\alpha}^R)=0.
$$

 Therefore, by the Lagrange multipliers rule, there exists $\Lambda_{\alpha,k}^D \in \mathbb{R}$ (resp. $\Lambda_{\alpha,k}^N \in \mathbb{R}, \Lambda_{\alpha,k}^R \in \mathbb{R})$ such that
$$
\mathcal{J}'(u_{\alpha,k}^D) = \Lambda_{\alpha,k}^D \mathcal{I}'(u_{\alpha,k}^D) \quad \text{ weakly in } \Omega,
$$
$$
(\text{resp.}\quad \mathcal{J}_0'(u_{\alpha,k}^N) = \Lambda_{\alpha,k}^N \mathcal{I}'(u_{\alpha,k}^N),\qquad  \mathcal{J}_{\beta}'(u_{\alpha,k}^R) = \Lambda_{\alpha,k}^R \mathcal{I}'(u_{\alpha,k}^R)),
$$
that is, $\{\Lambda^D_{k,\alpha}\}_{k\in\N}$, $\{\Lambda^N_{k,\alpha}\}_{k\in\N}$ and $\{\Lambda^R_{k,\alpha}\}_{k\in\N}$ are eigenvalues of \ref{eq.d}, \ref{eq.n} and \ref{eq.r}, respectively, satisfying that
$$
\mathcal{J}(u_{k,\alpha}^D)=C_{k,\alpha}^D, \qquad \mathcal{J}_0(u_{k,\alpha}^N)=C_{k,\alpha}^N, \qquad  \mathcal{J}_\beta (u_{k,\alpha}^R)=C_{k,\alpha}^R.
$$
The proof is concluded.
\end{proof}


Finally we provide for a proof of the comparison result Theorem \ref{compara}.

\begin{proof}[Proof of Theorem \ref{compara}]
Let $u_{k,\alpha}^D$ be the eigenfunctions corresponding to $\lambda_{k,\alpha}^D$. By definition of $\lambda_{k,\alpha}^D$ and property \eqref{G1} we get
$$
\lambda_{k,\alpha}^D= \frac{\langle (-\Delta_g)^s u_{k,\alpha}^D,u_{k,\alpha}^D \rangle}{\int_{\Omega} g(|u_{k,\alpha}^D|) |u_{k,\alpha}^D|\, dx} \leq \frac{p^+}{p^-} \frac{\Phi_{s,G,\R^n}(u_{k,\alpha}^D)}{\Phi_{G,\Omega}(u_{k,\alpha}^D)}= \frac{p^+}{p^-} \frac{\alpha}{c_{k,\alpha}^D}.
$$
In a similar way, if we consider the sequence $v_{k,\alpha}^D$ of eigenfunctions corresponding to $\Lambda_{k,\alpha}^D$, by using \eqref{G1} we get
$$
\Lambda_{k,\alpha}= \frac{\langle (-\Delta_g)^s v_{k,\alpha}^D,v_{k,\alpha}^D \rangle}{\int_{\Omega} g(|v_{k,\alpha}^D|) |v_{k,\alpha}^D|\, dx} \leq \frac{p^+}{p^-} \frac{\Phi_{s,G,\R^n}(v_{k,\alpha}^D)}{\Phi_{G,\Omega}(v_{k,\alpha}^D)}= \frac{p^+}{p^-} \frac{C_{k,\alpha}^D}{\alpha}.
$$
The lower bounds follow analogously.
\end{proof}

Finally, we provide for a proof of Theorem \ref{dirichlet.diverge}

\begin{proof}[Proof of Theorem \ref{dirichlet.diverge}]
{\bf Step 1.}
In \cite[Corollary 2.10]{FBPLS} it is proved that $W^{s,G}_0(\Omega)\subset W^{s,p^-}_0(\Omega)$, therefore
$$
[u]_{W^{s,p^-}(\R^n)} \leq C [u]_{W^{s,G}(\R^n)}.
$$
Given $u\in M_\alpha^D$ we get
$$
\mathcal{J}(u) = \Phi_{s,G,\R^n}(u) \geq \xi^-([u]_{W^{s,G}(\R^n)}) \geq C \xi^-([u]_{W^{s,p^-}(\R^n)})
$$
where we have used Lemma \ref{ineq1}. Then
\begin{equation} \label{eq11}
C[u]_{W^{s,p^-}(\R^n)}^\beta \leq \mathcal{J}(u)
\end{equation}
for some $\beta=\beta(p^\pm)$.

{\bf Step 2.} Define
$$
\mathcal{M}_\delta^D:=\{u\in W^{s,p^-}_0(\Omega)\colon [u]^\beta_{L^{p^-}(\Omega)} \leq \delta\}
$$
where $\beta$ is the same of step 1.

As in \cite[Lemma 2.7]{FBPLS} it can be seen that $L^G(\Omega)\subset L^{p^-}(\Omega)$. Then $\|u\|_{L^{p^-}(\Omega)} \leq C \|u\|_{L^G(\Omega)}$.

Given $u\in M_\alpha^D$,
$$
\alpha = \mathcal{I}(u) \geq \xi^-(\|u\|_{L^G(\Omega)}) \geq C\xi^-(\|u\|_{L^{p^-}(\Omega)}).
$$
Therefore, there exists some $\delta=\delta(\alpha,p^\pm)$ such that $\|u\|^\beta_{L^{p^-}(\Omega)}\leq \delta$ and then
\begin{equation} \label{eq12}
M_\alpha^D \subset \mathcal{M}_\delta^D.
\end{equation}

{\bf Step 3.}
By \eqref{eq11} and \eqref{eq12}
\begin{align} \label{eq13}
\begin{split}
C_{k,\alpha}^D &= \inf_{h\in \Gamma(S^{k-1},M_\alpha^D)} \sup_{w \in S^{k-1}} \mathcal{J}(h(w))\\
&\geq
 C \inf_{h\in \Gamma(S^{k-1},\mathcal{M}_\delta^D)} \sup_{w \in S^{k-1}}  [h(w)]_{W^{s,p^-}}^\beta := C(\mu_k^D)^\beta
\end{split}
\end{align}
where $\mu_k^D$ is the minimax eigenvalue of the $p^--$Laplacian with Dirichlet boundary conditions obtained in \cite[Theorem 4.1]{LS} (observe that since the $p^--$Laplacian is a homogeneous operator, in fact the same eigenvalue is obtained for any $\delta$).

{\bf Step 4.}
From Theorem \ref{compara} and \eqref{eq13} we get
$$
C\frac{p^-}{\alpha p^+} (\mu_k^D)^\beta \leq \frac{p^-}{\alpha p^+} C_{k,\alpha}^D \leq \Lambda_{k,\alpha}^D.
$$
Since $\mu_k^D \to \infty$ as $k\to \infty$ we obtain that $C_{k,\alpha}^D, \Lambda_{k,\alpha}^D\to\infty$ as $k\to\infty$.
\end{proof}

\subsection*{Acknowledgements.}The third author  was partially supported by Consejo Nacional de Investigaciones Cient\'{i}ficas y T\'{e}cnicas (CONICET-Argentina).

\end{document}